\documentclass{article}
\usepackage{amsmath,amssymb}	
\usepackage{graphicx}

\newtheorem{thm}{Theorem}
\newtheorem{lem}{Lemma}
\newtheorem{prop}{Proposition}
\newtheorem{cor}{Corollary}

\newtheorem{rem}{Remark}

\newcommand{\Dt}{D_{\tau}}

\title{Universal central extensions of linear groups over rings of
non-commutative Laurent polynomials, associated $K_1$-groups and $K_2$-groups}
\author{Ryusuke Sugawara}
\date{}

\begin{document}
\maketitle

\section{Introduction}

Many researchers have studied the structure of general linear groups and their elementary subgroups
over fields $F$ or commutative rings $R$. They also have analyzed associated lower $K$-groups, for example \cite{m} and \cite{s}. Needless to say, general linear groups are important objects and have many applications in various areas of mathematics, but they particularly have much to do with Lie theory; Lie groups, Lie algebras and their representations. Lower $K$-groups also play an important role as a certain invariant.

In this paper, we treat some rings $D_{\tau}$ of non-commutative Laurent polynomials over division rings $D$ (cf. Section 2). Here $\tau$ is an automorphism of $D$. We note that the ring we use generalizes the one which is studied in \cite{m} and \cite{s}. When $D=F$ and $\tau=id$, our discussion is just a subject of loop groups which are applied in the theory of affine Kac-Moody Lie algebras, and this is surveyed in \cite{jm} for example. On the other hand, the corresponding linear group was studied in the case when $D$ is the field of formal power series and $\tau$ is not trivial (cf. \cite{ms},\cite{hs}), which is deeply related to the theory of extended affine Lie algebras (cf. \cite{aabgp},\cite{mt},\cite{yy}). 

Our main object in this study is the following exact sequence. 
\[1\to K_2(n,D_{\tau})\to St(n,D_{\tau})\xrightarrow{\phi} GL(n,D_{\tau})\to K_1(n,D_{\tau})\to 1\]
We reveal the structure of groups in the above sequence.
We first describe an existence of a Tits system in the elementary subgroup $E(n,D_{\tau})$ of the general linear group $GL(n,D_{\tau})$ and the associated Steinberg group $St(n,D_{\tau})$ in Section 2 and Section 3 respectively. Using this fact, we show the above homomorphism $\phi$ is a central extension of $E(n,D_{\tau})$, that is, we confirm $\mathrm{Ker}\phi:=K_2(n,D_{\tau})$ is a central subgroup of the Steinberg group $St(n,\Dt)$ in Section 3. It is proved in Section 5 that $\phi$ is actually universal when the center $Z(D)$ of $D$ has  at least five elements. Meanwhile, we discuss the structure of the associated $K_1$-group and $K_2$-group in Section 4 and 6. In particular we check that $K_2$-groups are generated by certain products of Steinberg symbols. While we deal with linear groups of rank two for simplicity of discussion, we need to raise a rank of these groups for the sake of the general theory which is possible by \cite{mj}.

\begin{center}
\textbf{Acknowledgements}
\end{center}
The author would like to thank his advisor, Professor Jun Morita, for
his continuous guidance and encouragement.\vspace{3mm}

\section{Linear groups over rings of non-commutative Laurent polynomials}

Let $D$ be a division ring, and we fix an automorphism $\tau\in\mathrm{Aut}(D)$.
In the following, we denote by $D_{\tau}=D[t,t^{-1}]$ the ring of Laurent polynomials generated by
$D$ and an indeterminate $t$, whose ring structure is given by $tat^{-1}=\tau(a)$ for all $a\in D$.
Let $M(n,\Dt)$ be the ring of $n\times n$ matrices whose entries are in $\Dt$, and
we define the general linear group $GL(n,\Dt)$ as the matrix group of all invertible matrices,
that is, $GL(n,\Dt)=M(n,\Dt)^{\times}$.\vspace{3mm}

Let $\Delta=\{\epsilon_i-\epsilon_j\ |\ 1\leq i\neq j\leq n\}$ be a root system of type $\mathrm{A}_{n-1}$, where $\{ \epsilon_i\}_{1\leq i\leq n}$ is an orthonormal basis with respect to an inner product $(\cdot,\cdot)$ defined by $(\epsilon_i,\epsilon_j)=\delta_{ij}$. We see that any root in $\Delta$ is expressed as
\[ \epsilon_i-\epsilon_j=(\epsilon_i-\epsilon_{i+1})+(\epsilon_{i+1}-\epsilon_{i+2})+\cdots+(\epsilon_{j-1}-\epsilon_j)\]
if $i<j$ and as its minus version if $i>j$, so we put $\Pi=\{ \alpha_i:=\epsilon_i-\epsilon_{i+1}\ |\ 1\leq i\leq n-1\}$ and call it a simple system of $\Delta$. We call $\Delta^+:=\text{Span}_{\mathbb{Z}_{\geq 0}}\Pi\cap\Delta$ a set of positive roots and $\Delta^-:=-\Delta^+$ a set of negative roots. Also, let $\Delta_a:=\Delta\times\mathbb{Z}$ be an (abstract) affine root system of type $\mathrm{A}_{n-1}^{(1)}$ and $\Pi_a=\{ \dot{\alpha_i}:=(\alpha_i,0)|1\leq i\leq n-1\}\cup\{\dot{\alpha_0}:=(-\theta,1)\}$ be a simple system of $\Delta_a$, where $\theta=\alpha_1+\dots+\alpha_{n-1}$ is a highest root in $\Delta$. \vspace{3mm}

For $\beta=\epsilon_i-\epsilon_j\in\Delta$, $f\in \Dt$ we define
	\[ e_{\beta}(f)=I+fE_{ij}\]
where $I$ is the identity matrix and $E_{ij}$ is the matrix unit.
For $\dot{\beta}:=(\beta, m)\in\Delta_a$, $f\in D$ and $s,u\in D^{\times}=D\setminus\{ 0\}$ we put 
\begin{align*}
	& x_{\dot{\beta}}(f)=\begin{cases}
						e_{\beta}(ft^m) & (\beta\in\Delta^+), \\
						e_{\beta}(t^mf) & (\beta\in\Delta^-),
					 \end{cases}\\
	&w_{\dot{\beta}}(u)=x_{\dot{\beta}}(u)x_{-\dot{\beta}}(-u^{-1})x_{\dot{\beta}}(u),\\
	&h_{\dot{\beta}}(s)=w_{\dot{\beta}}(s)w_{\beta}(-1),
\end{align*}
then we can easily see $x_{\dot{\beta}}(f)^{-1}=x_{\dot{\beta}}(-f)$, $w_{\dot{\beta}}(u)^{-1}=w_{\dot{\beta}}(-u)$.
The elementary subgroup
$E(n,\Dt)$ is defined to be the subgroup of $GL(n,\Dt)$ generated by $x_{\dot{\beta}}(f)$ for all
$\dot{\beta}\in\Delta_a$ and $f\in D$.\vspace{3mm}

In a standard way,  the Weyl group $W$ of $\Delta$ is generated by all reflections $\sigma_{\beta}$ for $\beta\in\Delta$, 
and the Weyl group $W_a$ of $\Delta_a$ (the affine Weyl group of $\Delta$) is generated by all reflections
$\sigma_{\dot{\beta}}$ for $\dot{\beta}\in\Delta_a$, where the action of $\sigma_{\dot{\beta}}$ is defined as
\[\sigma_{\dot{\beta}}(\dot{\gamma})=(\sigma_{\beta}(\gamma), n-\langle \gamma, \beta \rangle m)\]
for $\dot{\beta}=(\beta, m),\ \dot{\gamma}=(\gamma, n)\in\Delta_a$ and 
$\langle \gamma,\beta\rangle=2(\gamma,\beta)/(\beta,\beta)$.
Define $\xi_{\dot{\beta}}:=\sigma_{\dot{\beta}}\sigma_{\beta}$ for each $\dot{\beta}\in\Delta_a$,
then let $T_a$ be the subgroup of $W_a$ generated by $\xi_{\dot{\beta}}$ for all $\dot{\beta}\in\Delta_a$.
In order to use later, we will characterize the group $W_a$ as follow (cf. \cite{jm} Lemma 1.1, Proposition 1.2).\vspace{3mm}

{\lem
\begin{align*}
(1)\ &\text{Let $\dot{\beta}=(\beta,m)$ and $\dot{\gamma}=(\gamma,n)$. Then }
\xi_{\dot{\beta}}(\dot{\gamma})=(\gamma, n+\langle\gamma,\beta\rangle m). \\
(2)\ &\text{Let $\alpha\in \Pi$. Then $T_a$ is a free abelian group generated by $\xi_{\dot{\gamma}}$ for all $\dot{\gamma}=(\alpha,1)$.}\\
(3)\ &\text{$\sigma_{\dot{\beta}}$ normalizes $T_a$ for $\dot{\beta}\in\Delta_a$. $\square$} 
\end{align*}
}

{\prop $W_a\cong T_aW\cong T_a\rtimes W$. $\square$} \vspace{5mm}

As a subgroup of $E(n,\Dt)$ we put
\begin{align*}
	&U_{\dot{\beta}}=\{ x_{\dot{\beta}}(f)\ |\ f\in D\},\\
	&U^{\pm}=\langle U_{\dot{\beta}}\ |\ \dot{\beta}\in\Delta_a^{\pm}\rangle,\\
	&N=\langle w_{\dot{\beta}}(u)\ |\ \dot{\beta}\in\Delta_a,\ u\in D^{\times}\rangle,\\
	&T=\langle h_{\dot{\beta}}(u)\ |\ \dot{\beta}\in\Delta_a,\ u\in D^{\times}\rangle,
\end{align*}
where $\Delta_a^+=(\Delta^+\times\mathbb{Z}_{\geq 0})\cup(\Delta^-\times\mathbb{Z}_{>0})$ and
$\Delta_a^-=(\Delta^+\times\mathbb{Z}_{< 0})\cup(\Delta^-\times\mathbb{Z}_{\leq 0})$.
If $h\in T$ is expressed as $h=\mathrm{diag}(u_1,u_2,\dots,u_n)$, the diagonal matrix with $u_i\in\Dt^{\times}=D^{\times}\cdot \{ t^l|l\in\mathbb{Z}\}$, then we define
	\[ \mathrm{deg}_i(h)=\mathrm{deg}(u_i)=m_i \]
for $i=1,2,\dots,n$, where $u_i=s_it^{m_i}$ with $s_i\in D^{\times}$ and $m_i\in\mathbb{Z}$. Then we set
\begin{align*}
	&T_0=\langle h\ |\ h\in T,\ \mathrm{deg}_i(h)=0\ \mathrm{for\ all}\ i=1,2,\dots,n\rangle,\\
	&B^{\pm}=\langle U^{\pm},\ T_0\rangle,\\
	&S=\{ w_{\dot{\beta}}(1)\ \mathrm{mod}\ T_0\ |\ \dot{\beta}\in\Delta_a\}.
\end{align*}
In particular $S$ is identified with the set $\{ w_{\dot{\alpha}}\ |\ \dot{\alpha}\in\Pi_a\}$.
The main result in this section is the following theorem.\vspace{3mm}

{\thm{Notation is as above. $(E(n,\Dt),B^{\pm},N,S)$ is a Tits system with the corresponding affine Weyl group $W_a$.}}\vspace{5mm}

Before proving this theorem  we give several relations between $e$, $w$ and
$h$.
\begin{align*}
&(R1)\ e_{\beta}(f)e_{\beta}(g)=e_{\beta}(f+g),\\
&(R2)\ [e_{\beta}(f),e_{\gamma}(g)]=
	\begin{cases}
		e_{\beta+\gamma}(fg) \hspace{5mm}&\text{if}\ \beta+\gamma\in\Delta,\ j=k,\\
		e_{\beta+\gamma}(-gf) \hspace{5mm}&\text{if }\ \beta+\gamma\in\Delta,\ i=l,\\
		1 \hspace{10mm} &\text{otherwise},
	\end{cases}
\end{align*}
where  $\beta=\epsilon_i-\epsilon_j$, $\gamma=\epsilon_k-\epsilon_l$, $i\neq l,j\neq k$ and $f,g\in D_{\tau}$.
These are fundamental relations in $E(n,\Dt)$ for example. Then we obtain;
\begin{align*}
&(R3)\ w_{\beta}(u)e_{\gamma}(f)w_{\beta}(u)^{-1}=
	\begin{cases}
		e_{\gamma}(f)\hspace{5mm}&\text{if} \hspace{2mm}(\beta,\gamma)=0,\\
		e_{\mp\beta}(u^{\mp 1}fu^{\mp 1})\hspace{5mm}&\text{if} \hspace{2mm}\gamma=\pm\beta,\\
		e_{\sigma_{\beta}(\gamma)}(-u^{-1}f)\hspace{5mm}&\text{if} \hspace{2mm}\beta\pm\gamma\neq 0\ \text{and}\ i=k,\\
		e_{\sigma_{\beta}(\gamma)}(-fu)\hspace{5mm}&\text{if} \hspace{2mm}\beta\pm\gamma\neq 0\ \text{and}\ i=l,\\
		e_{\sigma_{\beta}(\gamma)}(uf)\hspace{5mm}&\text{if} \hspace{2mm}\beta\pm\gamma\neq 0\ \text{and}\ j=k,\\
		e_{\sigma_{\beta}(\gamma)}(fu^{-1})\hspace{5mm}&\text{if} \hspace{2mm}\beta\pm\gamma\neq 0\ \text{and}\ j=l,\\
	\end{cases}\\
&(R4)\ h_{\beta}(u)e_{\gamma}(f)h_{\beta}(u)^{-1}=
	\begin{cases}
		e_{\gamma}(f)\hspace{5mm}&\text{if} \hspace{2mm}(\beta,\gamma)=0,\\
		e_{\pm\beta}(-u^{\pm 1}fu^{\pm 1})\hspace{5mm}&\text{if} \hspace{2mm}\gamma=\pm\beta,\\
		e_{\gamma}(uf)\hspace{5mm}&\text{if} \hspace{2mm}\beta\pm\gamma\neq 0\ \text{and}\ i=k,\\
		e_{\gamma}(fu^{-1})\hspace{5mm}&\text{if} \hspace{2mm}\beta\pm\gamma\neq 0\ \text{and}\ i=l,\\
		e_{\gamma}(u^{-1}f)\hspace{5mm}&\text{if} \hspace{2mm}\beta\pm\gamma\neq 0\ \text{and}\ j=k,\\
		e_{\gamma}(fu)\hspace{5mm}&\text{if} \hspace{2mm}\beta\pm\gamma\neq 0\ \text{and}\ j=l,\\
	\end{cases}\\
&(R5)\ w_{\beta}(u)w_{\gamma}(s)w_{\beta}(u)^{-1}=
	\begin{cases}
		w_{\gamma}(s)\hspace{5mm}&\text{if} \hspace{2mm}(\beta,\gamma)=0,\\
		w_{\mp\beta}(-u^{\mp 1}su^{\mp 1})\hspace{5mm}&\text{if} \hspace{2mm}\gamma=\pm\beta,\\
		w_{\sigma_{\beta}(\gamma)}(-u^{-1}s)\hspace{5mm}&\text{if} \hspace{2mm}\beta\pm\gamma\neq 0\ \text{and}\ i=k,\\
		w_{\sigma_{\beta}(\gamma)}(-su)\hspace{5mm}&\text{if} \hspace{2mm}\beta\pm\gamma\neq 0\ \text{and}\ i=l,\\
		w_{\sigma_{\beta}(\gamma)}(us)\hspace{5mm}&\text{if} \hspace{2mm}\beta\pm\gamma\neq 0\ \text{and}\ j=k,\\
		w_{\sigma_{\beta}(\gamma)}(su^{-1})\hspace{5mm}&\text{if} \hspace{2mm}\beta\pm\gamma\neq 0\ \text{and}\ j=l,\\
	\end{cases}\\
&(R6)\ w_{\beta}(u)h_{\gamma}(s)w_{\beta}(u)^{-1}=
	\begin{cases}
		h_{\gamma}(s)\hspace{5mm}&\text{if} \hspace{2mm}(\beta,\gamma)=0,\\
		h_{\mp\beta}(u^{\mp 1}su^{\mp 1})h_{\mp\beta}(u^{\pm 2})\hspace{5mm}&\text{if} \hspace{2mm}\gamma=\pm\beta,\\
		h_{\sigma_{\beta}(\gamma)}(u^{-1}s)h_{\sigma_{\beta}(\gamma)}(u)\hspace{5mm}&\text{if} \hspace{2mm}\beta\pm\gamma\neq 0\ \text{and}\ i=k,\\
		h_{\sigma_{\beta}(\gamma)}(su)h_{\sigma_{\beta}(\gamma)}(u^{-1})\hspace{5mm}&\text{if} \hspace{2mm}\beta\pm\gamma\neq 0\ \text{and}\ i=l,\\
		h_{\sigma_{\beta}(\gamma)}(us)h_{\sigma_{\beta}(\gamma)}(u^{-1})\hspace{5mm}&\text{if} \hspace{2mm}\beta\pm\gamma\neq 0\ \text{and}\ j=k,\\
		h_{\sigma_{\beta}(\gamma)}(su^{-1})h_{\sigma_{\beta}(\gamma)}(u)\hspace{5mm}&\text{if} \hspace{2mm}\beta\pm\gamma\neq 0\ \text{and}\ j=l,\\
	\end{cases}
\end{align*}
for all $\beta,\gamma\in\Delta$, $f,g\in D_{\tau}$ and $s,u\in D_{\tau}^{\times}$, where $\beta=\epsilon_i-\epsilon_j$ and $\gamma=\epsilon_k-\epsilon_l$.

{\lem\begin{align*}
		&(1)\ B^{\pm}=U^{\pm}\rtimes T_0,\\
		&(2)\ T_0\lhd N\ and\ T\lhd N,\\
		&(3)\ B^{\pm}\cap N=T_0,\\
		&(4)\ N/T_0\cong W_a.
			\end{align*}}
\textbf{Proof:}\ From the definition, we know $B^{\pm}=\langle U^{\pm},T_0\rangle$,
and we also see $U^{\pm}\cap T_0=\{ I\}$ easily
if we consider the degree of an element in $U^{\pm}$.
In addition $T_0$ normalizes $U^{\pm}$ from (R4), hence (1) holds. (2) has already been proven by (R6).
We prove (4) along with \cite{jm}. If we put
$N_0=\langle w_{\alpha}(u)| \alpha\in\Pi,u\in D^{\times}\rangle$ then we have $N=TN_0$ by the definition of $T$, and
we know $T_0\lhd T$ and $\ T_0\lhd N$. Therefore
\[N/T_0\cong (T/T_0)(N_0/T_0)\cong (T/T_0)W.\]
Here we set $\dot{\gamma}=(\alpha,1)$ for $\alpha\in\Pi$. Then the map $T_a\to T/T_0$ defined by 
$\xi_{\dot{\gamma}}\mapsto h_{\alpha}(t)$ is an isomorphism,
hence $N/T_0\cong W_a$ from Proposition 1. By $U^{\pm}\cap T_0=\{ I\}$ and (4),
we can see that (3) holds. $\square$\vspace{3mm}

To discuss these subgroups more explicitly we put 
\[ U'_{\pm\dot{\alpha}}=\langle x_{\pm\dot{\alpha}}(g)U_{\dot{\beta}}x_{\pm\dot{\alpha}}(g)^{-1}\ |\ g\in D,\dot{\beta}\in\Delta^{\pm}_a\setminus\{ \pm\dot{\alpha}\}\rangle\]
for each $\dot{\alpha}\in\Pi_a$.\vspace{3mm}

{\prop Let $\dot{\alpha}\in\Pi_a$. Then;
\begin{align*}
(1)\ &w_{\pm\dot{\alpha}}(u)U'_{\pm\dot{\alpha}}w_{\pm\dot{\alpha}}(u)^{-1}=U'_{\pm\dot{\alpha}}\ \text{for all $u\in D^{\times}$},\\
(2)\ &U^{\pm}=U'_{\pm\dot{\alpha}}\rtimes U_{\pm\dot{\alpha}}.
\end{align*}} 
\textbf{Proof:}
(1) It suffices to show that 
\[w_{\dot{\alpha}}(u)x_{\dot{\alpha}}(g)x_{\dot{\beta}}(f')x_{\dot{\alpha}}(g)^{-1}w_{\dot{\alpha}}(u)^{-1}\in U_{\dot{\alpha}}'\]
for all $\dot{\beta}\in\Delta_a^+\setminus\{\dot{\alpha}\}$. For $\dot{\beta}=(\beta,m)$ we denote its negative root in a first entry by $\dot{-\beta}$, that is, $\dot{-\beta}=(-\beta,m)$.
There are two cases:\vspace{3mm}

Case 1. If $\dot{\beta}=(\beta, m)\in\Delta_a^+\setminus\{ \dot{\alpha}\}$ and $\beta\neq -\alpha$, then we see
\begin{align*}
w_{\dot{\alpha}}(u)x_{\dot{\alpha}}(g)x_{\dot{\beta}}(f')x_{\dot{\alpha}}(g)^{-1}w_{\dot{\alpha}}(u)^{-1}&=
	\begin{cases}
		w_{\dot{\alpha}}(u)x_{\dot{\alpha}+\dot{\beta}}(uf')x_{\dot{\beta}}(f')w_{\dot{\alpha}}(u)^{-1}\\
		w_{\dot{\alpha}}(u)x_{\dot{\alpha}+\dot{\beta}}(-f'u)x_{\dot{\beta}}(f')w_{\dot{\alpha}}(u)^{-1}\\
		w_{\dot{\alpha}}(u)x_{\dot{\beta}}(f')w_{\dot{\alpha}}(u)^{-1}
	\end{cases}\\&=
	\begin{cases}
		x_{\sigma_{\dot{\alpha}}(\dot{\alpha}+\dot{\beta})}(*)x_{\sigma_{\dot{\alpha}}(\dot{\beta})}(*)\\
		x_{\sigma_{\dot{\alpha}}(\dot{\alpha}+\dot{\beta})}(*)x_{\sigma_{\dot{\alpha}}(\dot{\beta})}(*)\\
		x_{\sigma_{\dot{\alpha}}(\dot{\beta})}(*)
	\end{cases}
\end{align*}
by the previous relations $(R2)$ and $(R3)$, but these all belong to $U'_{\dot{\alpha}}$ because of the fact $m\in\mathbb{Z}_{>0}$ in any cases.\vspace{3mm}

Case 2. If $\dot{\beta}=(\beta, m)\in\Delta_a^+\setminus\{ \dot{\alpha}\}$ and $\beta= -\alpha$, then
\begin{align*}
&w_{\dot{\alpha}}(u)x_{\dot{\alpha}}(g)x_{\dot{\beta}}(f')x_{\dot{\alpha}}(g)^{-1}w_{\dot{\alpha}}(u)^{-1}\\
&=w_{\dot{\alpha}}(u)x_{\dot{\alpha}}(g)w_{\dot{\alpha}}(u)^{-1}w_{\dot{\alpha}}(u)x_{\dot{\beta}}(f')w_{\dot{\alpha}}(u)^{-1}w_{\dot{\alpha}}(u)x_{\dot{\alpha}}(g)^{-1}w_{\dot{\alpha}}(u)^{-1}\\
&=x_{\dot{-\alpha}}(g')x_{\dot{-\beta}}(f'')x_{\dot{-\alpha}}(g')^{-1}\\
&=x_{\dot{\alpha}}(g'^{-1})w_{\dot{\alpha}}(-g'^{-1})x_{\dot{\alpha}}(g'^{-1})x_{\dot{-\beta}}(f'')x_{\dot{\alpha}}(g'^{-1})^{-1}w_{\dot{\alpha}}(-g'^{-1})^{-1}x_{\dot{\alpha}}(g'^{-1})^{-1}\\
&=x_{\dot{\alpha}}(g'^{-1})x_{\dot{\beta}}(f''')x_{\dot{\alpha}}(g'^{-1})^{-1}\in U_{\dot{\alpha}}'.
\end{align*}
(2) From the definition, we see $U'_{\pm\dot{\alpha}}<U^{\pm}$ and $U_{\dot{\pm\alpha}}<U^{\pm}$, and that
$U_{\pm\dot{\alpha}}$ normalizes $U'_{\dot{\pm\alpha}}$.
Suppose $U'_{\pm\dot{\alpha}}\cap U_{\pm\dot{\alpha}}\neq \{ I\}$.
Then, for each $x_{\pm\dot{\alpha}}(f)\in U_{\pm\dot{\alpha}}$, there exist $x'_i\in U'_{\pm\dot{\alpha}}$
such that $x_{\pm\dot{\alpha}}(f)=x'_1x'_2\dots x'_r$.
By (1), we have
\[U^{\mp}\ni w_{\pm\dot{\alpha}}(u)x_{\pm\dot{\alpha}}(f)w_{\pm\dot{\alpha}}(u)^{-1}=w_{\pm\dot{\alpha}}(u)x'_1x'_2\dots x'_rw_{\pm\dot{\alpha}}(u)^{-1}\in U^{\pm}.\]
Thus $U'_{\pm\dot{\alpha}}\cap U_{\pm\dot{\alpha}}= \{ I\}$. $\square$\vspace{3mm}

We need to check the following two conditions for the sake of the proof; 
\begin{align*}
(T1)\ &sB^{\pm}w\subset B^{\pm}wB^{\pm}\cup B^{\pm}swB^{\pm}\ \text{for all $s\in S$ and $w\in W_a$},\\
(T2)\ &sB^{\pm}s\not\subset B^{\pm}\ \text{for each $s\in S$}.
\end{align*}

\textbf{Proof of Theorem 1}: We check the ``axiom of Tits systems''. We have  already see that $E(n,\Dt)$ is a group,
$B^{\pm}$ and $N$ are subgroups of $E(n,\Dt)$ and $S$ is a subset of $N/(B^{\pm}\cap N)$.
In addition $E(n,\Dt)=\langle B^{\pm},N\rangle$ and $B^{\pm}\cap N\lhd N$.
Moreover $N/(B^{\pm}\cap N)=W_a=\langle S\rangle$, and all elements in $S$ are of order two (modulo $T_0$).
Finally we prove the above two conditions (T1) and (T2), but it is enough to check only the case of $B^+$. In the following,
we write $B=B^+$.\vspace{3mm}

Define the length, called $l(w)$, of an element $w$ in $W_a$ as usual. For $w\in W_a$ and $s=w_{\dot{\alpha}}(1)\in S$,
if $l(w)<l(sw)$ then we have $w(\dot{\alpha})\in\Delta_a^+$, hence
\begin{align*}
wBs&=wU_{\dot{\alpha}}U'_{\dot{\alpha}}T_0s\\
		&=wU_{\dot{\alpha}}w^{-1}wsU'_{\dot{\alpha}}sT_0s\\
		&=U_{w(\dot{\alpha})}wsU'_{\dot{\alpha}}T_0\\
		&\subset BwsB.
\end{align*}
Also, if $l(ws)<l(w)$ then $w(\dot{\alpha})\in\Delta_a^-$, and we have $l(w')<l(w's)$ for $w'=ws$. Thus
\begin{align*}
wBs&=w'sBs\\
		&\subset w'(B\cup BsB)\\
		&=w'B\cup w'BsB\\
		&\subset Bw'B\cup Bw'sBB\\
		&=BwsB\cup BwB,
\end{align*}
so (T1) holds. On the other hand, we can easily see that (T2) holds by direct calculation. $\square$\vspace{3mm}

{\cor{Notation is as above. Then\vspace{3mm}\\
$(1)$ $E(n,\Dt)$ has a Bruhat decomposition: 
\begin{center}$E(n,\Dt)=\bigcup_{w\in W_a}B^{\pm}wB^{\pm}=U^{\pm}NU^{\pm}$,\end{center}
$(2)$ $E(n,\Dt)$ has a Birkhoff decomposition:
\begin{center}$E(n,\Dt)=\bigcup_{w\in W_a}B^{\mp}wB^{\pm}=U^{\mp}NU^{\pm}$,\end{center}
(3) $E(n,\Dt)$ has a Gauss decomposition:
\begin{center}$E(n,\Dt)=\bigcup_{X\in U^{\pm}}XB^{\mp}B^{\pm}X^{-1}=U^{\pm}B^{\mp}U^{\pm}$. $\square$\end{center}}} 

\section{Steinberg groups}
Let $St(n,\Dt)$ be the Steinberg group, which is defined by the generators $\hat{x}_{ij}(f)$ for all $f\in\Dt$,
$1\leq i\neq j\leq n$ and the defining relations
\begin{align*}
&(ST1)\ \hat{x}_{ij}(f)\hat{x}_{ij}(g)=\hat{x}_{ij}(f+g)\\
&(ST2)\ [\hat{x}_{ij}(f),\hat{x}_{kl}(g)]=
		\begin{cases}
			\hat{x}_{il}(fg) &\text{if   } j=k,\\
			\hat{x}_{kj}(-gf) &\text{if   } i=l, \\
			1 &\text{otherwise},
		\end{cases}
\end{align*}
where $f,g\in\Dt$, $1\leq i\neq j\leq 2$ and $1\leq k\neq l\leq n$ with$(i,j)\neq (l,k)$. Then, we see $\hat{x}_{ij}(f)^{-1}=\hat{x}_{ij}(-f)$
from $(ST1)$. When $n=2$, we use
\[(ST2)'\ \hat{w}_{ij}(u)\hat{x}_{ij}(v)\hat{w}_{ij}(-u)=\hat{x}_{ji}(-u^{-1}vu^{-1})\]
instead of $(ST2)$, where $\hat{w}_{ij}(u)=\hat{x}_{ij}(u)\hat{x}_{ji}(-u^{-1})\hat{x}_{ij}(u)$.\vspace{3mm}

We now consider a natural homomorphism $\phi$ from $St(n,\Dt)$ onto $E(n,\Dt)$ defined by
$\phi(\hat{x}_{ij}(f))=I+E_{ij}$ for all $f\in\Dt$. Then, we claim that $\phi$ is a universal central extension of $E(n,\Dt)$
under some conditions. The proof is formed in three steps: we first show that
$St(n,\Dt)$ has a Tits system, and we check that $\phi$ is a central extension. Finally,
we prove the universality of $\phi$. In this section we show $St(n,\Dt)$ has a Tits system.\vspace{3mm}

For $\dot{\beta}=(\beta,m)\in\Delta_a$ and $f\in D$, we put
\[\hat{x}_{\dot{\beta}}(f)=
				 \begin{cases}
					\hat{x}_{\beta}(ft^m) & (\beta\in\Delta^+), \\
					\hat{x}_{\beta}(t^mf) & (\beta\in\Delta^-),
				 \end{cases}\]
and we also put $\hat{w}_{\dot{\beta}}(u)=\hat{x}_{\dot{\beta}}(u)\hat{x}_{-\dot{\beta}}(-u^{-1})\hat{x}_{\dot{\beta}}(u)$
and $\hat{h}_{\dot{\beta}}(u)=\hat{w}_{\dot{\beta}}(u)\hat{w}_{\beta}(-1)$ for $u\in D^{\times}$ as in Section 2.
Here, we put, as a subgroup of $St(n,\Dt)$,
\begin{align*}
	&\hat{U}_{\dot{\beta}}=\{ \hat{x}_{\dot{\beta}}(f)\ |\ f\in D\},\\
	&\hat{U}^{\pm}=\langle \hat{U}_{\dot{\beta}}\ |\ \dot{\beta}\in\Delta_a^{\pm}\rangle,\\
	&\hat{N}=\langle \hat{w}_{\dot{\beta}}(u)\ |\ \dot{\beta}\in\Delta_a,\ u\in D^{\times}\rangle,\\
	&\hat{T}=\langle \hat{h}_{\dot{\beta}}(u)\ |\ \dot{\beta}\in\Delta_a,\ u\in D^{\times}\rangle.
\end{align*}
Then we set
\begin{align*}
	&\hat{T_0}=\langle h\ |\ h\in \hat{T},\ \mathrm{deg}_i(\phi(h))=0\ \mathrm{for\ all}\ i=1,\dots, n\rangle,\\
	&\hat{B}^{\pm}=\langle \hat{U}^{\pm},\ \hat{T_0}\rangle,\\
	&\hat{S}=\{ \hat{w}_{\dot{\beta}}(1)\ \mathrm{mod}\ \hat{T_0}\ |\ \dot{\beta}\in\Delta_a\}.
\end{align*}
In particular $\hat{S}$ is identified with the set $\{ \hat{w}_{\dot{\alpha}}\ |\ \dot{\alpha}\in\Pi_a\}$.\vspace{3mm}

We also need several relations between $\hat{x}_{\dot{\beta}}(f)$, $\hat{w}_{\dot{\beta}}(u)$ and $\hat{h}_{\dot{\beta}}(s)$
like $(R1)\sim (R6)$, but those actually are the same as in Section 2 except $(R6)$. We give $(R6)$ in Steinberg groups as follows (cf. \cite{s}):

\begin{align*}
(\hat{R}6)\ &\hat{w}_{\beta}(u)\hat{h}_{\gamma}(s)\hat{w}_{\beta}(u)^{-1}\\&=
	\begin{cases}
		\hat{h}_{\gamma}(s)\hspace{2mm}&\text{if} \hspace{2mm}(\beta,\gamma)=0,\\
		\hat{h}_{\mp\beta}(-u^{\mp 1}su^{\mp 1})\hat{h}_{\mp\beta}(-u^{\pm 2})^{-1} &\text{if} \hspace{2mm}\gamma=\pm\beta,\\
		\hat{h}_{\sigma_{\beta}(\gamma)}(-u^{-1}s)\hat{h}_{\sigma_{\beta}(\gamma)}(-u^{-1})^{-1} &\text{if} \hspace{2mm}\beta\pm\gamma\neq 0\ \text{and}\ i=k,\\
		\hat{h}_{\sigma_{\beta}(\gamma)}(-su)\hat{h}_{\sigma_{\beta}(\gamma)}(-u)^{-1} &\text{if} \hspace{2mm}\beta\pm\gamma\neq 0\ \text{and}\ i=l,\\
		\hat{h}_{\sigma_{\beta}(\gamma)}(us)\hat{h}_{\sigma_{\beta}(\gamma)}(u)^{-1} &\text{if} \hspace{2mm}\beta\pm\gamma\neq 0\ \text{and}\ j=k,\\
		\hat{h}_{\sigma_{\beta}(\gamma)}(su^{-1})\hat{h}_{\sigma_{\beta}(\gamma)}(u^{-1})^{-1} &\text{if} \hspace{2mm}\beta\pm\gamma\neq 0\ \text{and}\ j=l,\\
	\end{cases}
\end{align*}
where $\beta=\epsilon_i-\epsilon_j$, $\gamma=\epsilon_k-\epsilon_l\in\Delta$ and $u,s\in\Dt^{\times}$.
Then we have the following.\vspace{3mm}

{\lem\begin{align*}
		&(1)\ \hat{B}^{\pm}=\hat{U}^{\pm}\rtimes \hat{T_0},\\
		&(2)\ \hat{T_0}\lhd \hat{N}\ and\ \hat{T}\lhd \hat{N},\\
		&(3)\ \hat{B}^{\pm}\cap \hat{N}=\hat{T_0},\\
		&(4)\ \hat{N}/\hat{T_0}\cong W_a.\hspace{5mm}\square
			\end{align*}}

To discuss these subgroups more explicitly we put 
\[ \hat{U}'_{\pm\dot{\alpha}}=\langle \hat{x}_{\pm\dot{\alpha}}(g)\hat{U}_{\dot{\beta}}\hat{x}_{\pm\dot{\alpha}}(g)^{-1}\ |\ g\in D,\dot{\beta}\in\Delta^{\pm}_a\setminus\{ \pm\dot{\alpha}\}\rangle\]
for each $\dot{\alpha}\in\Pi_a$. Then the following proposition holds.\vspace{3mm}

{\prop\begin{align*}
	(1)\ &\hat{w}_{\pm\dot{\alpha}}(u)\hat{U}'_{\pm\dot{\alpha}}\hat{w}_{\pm\dot{\alpha}}(u)^{-1}=\hat{U}'_{\pm\dot{\alpha}}\ 		\text{for all $u\in D^{\times}$},\\
	(2)\ &\hat{U}^{\pm}=\hat{U}'_{\pm\dot{\alpha}}\hat{U}_{\pm\dot{\alpha}}=\hat{U}_{\pm\dot{\alpha}}\hat{U}'_{\pm\dot{\alpha}},\\
	(3)\ &\hat{s}\hat{B}^{\pm}w\subset \hat{B}^{\pm}w\hat{B}^{\pm}\cup \hat{B}^{\pm}\hat{s}w\hat{B}^{\pm}\ \text{for all $\hat{s}\in \hat{S}$ and $w\in W_a$},\\
	(4)\ &\hat{s}\hat{B}^{\pm}\hat{s}\not\subset \hat{B}^{\pm}\ \text{for each $\hat{s}\in \hat{S}$}.\hspace{5mm}\square
			\end{align*}} 

{\thm{Notation is as above. $(St(n,\Dt),\hat{B}^{\pm},\hat{N},\hat{S})$ is a Tits system with the corresponding affine Weyl group $W_a$.\hspace{5mm}$\square$}}\vspace{5mm}

{\cor{Notation is as above. Then\vspace{3mm}\\
$(1)$ $St(n,\Dt)$ has a Bruhat decomposition: 
\begin{center}$St(n,\Dt)=\bigcup_{w\in W_a}\hat{B}^{\pm}w\hat{B}^{\pm}=\hat{U}^{\pm}\hat{N}\hat{U}^{\pm}$,\end{center}
$(2)$ $St(n,\Dt)$ has a Birkhoff decomposition:
\begin{center}$St(n,\Dt)=\bigcup_{w\in W_a}\hat{B}^{\mp}w\hat{B}^{\pm}=\hat{U}^{\mp}\hat{N}\hat{U}^{\pm}$,\end{center}
(3) $St(n,\Dt)$ has a Gauss decomposition:
\begin{center}$St(n,\Dt)=\bigcup_{\hat{X}\in\hat{U}^{\pm}}\hat{X}\hat{B}^{\mp}\hat{B}^{\pm}\hat{X}^{-1}%
=\hat{U}^{\pm}\hat{B}^{\mp}\hat{U}^{\pm}$.\hspace{3mm}$\square$\end{center}}} 

\section{Central extensions}

We will show that a given natural homomorphism $\phi$ in Section 4 is a central extension, i.e.,
$\mathrm{Ker}\phi$ is a central subgroup of $St(n,\Dt)$. For this homomorphism $\phi$
we define our $K_2$-groups as $K_2(n,\Dt)=\mathrm{Ker}\phi$ (cf.\cite{m}). We first discuss a central extension of
$E(2,\Dt)$.\vspace{3mm}

Let $\tilde{E}(2,\Dt)$ be the group presented by the generators $\tilde{x}_{ij}(f)$ for all $f\in\Dt$, $1\leq i\neq j\leq 2$
and defining relations $(ST1)$, $(ST2)'$ and 
\begin{align*}
(ST3)\ \tilde{c}(u_1,v_1)\tilde{c}(u_2,v_2)\dotsm\tilde{c}(u_r,v_r)=1
\end{align*}
for all $u_i,v_i\in\Dt^{\times}$ such that $[u_1,v_1][u_2,v_2]\dotsm[u_r,v_r]=1$, where for $u,v\in\Dt^{\times}$ we put
\begin{align*}
&\tilde{w}_{ij}(u)=\tilde{x}_{ij}(u)\tilde{x}_{ji}(-u^{-1})\tilde{x}_{ij}(u),\\
&\tilde{h}_{ij}(u)=\tilde{w}_{ij}(u)\tilde{w}_{ij}(-1),\\
&\tilde{c}(u,v)=\tilde{h}_{12}(u)\tilde{h}_{12}(v)\tilde{h}_{12}(vu)^{-1},
\end{align*}
and where we change $\hat{x}$ and $\hat{w}$ in the defining relations to $\tilde{x}$ and $\tilde{w}$ respectively. That is, $\tilde{E}(2,\Dt)$ is the quotient group of $St(2,\Dt)$ by the corresponding elements $\hat{c}(u_1,v_1)\hat{c}(u_2,v_2)\dotsm\hat{c}(u_r,v_r)$.
Here, we note
\[\phi(\hat{c}(u,v))=\left(\begin{array}{cc}
	[u,v] & 0 \\
	0 & 1
	\end{array}\right),\]
and, in particular, we have $\hat{c}(u,v)\in\hat{T}_0$ and 
 \begin{align*}
\phi(\hat{c}(u_1,v_1)&\hat{c}(u_2,v_2)\dotsm\hat{c}(u_r,v_r))\\
&=
	\left(\begin{array}{cc}
	[u_1,v_1] & 0 \\
	0 & 1
	\end{array}\right)
	\left(\begin{array}{cc}
	[u_2,v_2] & 0 \\
	0 & 1
	\end{array}\right)\dots
	\left(\begin{array}{cc}
	[u_r,v_r] & 0 \\
	0 & 1
	\end{array}\right)\\
&=\left(\begin{array}{cc}
	[u_1,v_1][u_2,v_2]\dots [u_r,v_r] & 0 \\
	0 & 1
	\end{array}\right)\\
&=I
\end{align*}
if $[u_1,v_1][u_2,v_2]\dots [u_r,v_r]=1$.
Then we obtain the following theorem.\vspace{3mm}

{\thm{Notation is as above. We have $\tilde{E}(2,\Dt)\cong E(2,\Dt)$.}}\vspace{3mm}\\
\textbf{Proof of Theorem 3:} We will show this along with \cite{ms}. The homomorphism $\phi$ induces
two canonical homomorphism called $\hat{\phi}$ and $\tilde{\phi}$:
\begin{align*}
&\hat{\phi}: St(2,\Dt)\to \tilde{E}(2,\Dt),\\
&\tilde{\phi}: \tilde{E}(2,\Dt)\to E(2,\Dt),
\end{align*}
which are defined by 
\[ \hat{\phi}(\hat{x}_{ij}(f))=\tilde{x}_{ij}(f)\ \mathrm{and}\ \tilde{\phi}(\tilde{x}_{ij}(f))=I+fE_{ij}.\]
We use the same notation of subgroups of $\tilde{E}(2,\Dt)$ as in $St(2,\Dt)$ changing $\hat{}$ to $\tilde{}$.
Then, for any $\tilde{x}\in\mathrm{Ker}\tilde{\phi}$, we can choose $\tilde{y}\in\tilde{U}$ and $\tilde{z}\in\tilde{T_0}$
satisfying $\tilde{x}=\tilde{y}\tilde{z}$ because $St(2,\Dt)$ has a Bruhat decomposition and $\mathrm{Ker}\tilde{\phi}\subset \tilde{B}=\tilde{U}\tilde{T_0}$. We know
$\tilde{\phi}(\tilde{x})=\tilde{\phi}(\tilde{y})\tilde{\phi}(\tilde{z})=I$ and $U\cap T_0=\{ I\}$, so
we have $\tilde{\phi}(\tilde{y})=I$ and $\tilde{\phi}(\tilde{z})=I$, in particular, $\tilde{y}$ and $\tilde{z}$
belong to the kernel of $\tilde{\phi}$. Therefore we need to prove $\tilde{y}=1=\tilde{z}$.\vspace{3mm}

\textbf{STEP 1.} We claim $\tilde{z}=1$.

We know $\tilde{w}_{21}(u)=\tilde{w}_{12}(-u^{-1})$ from $(\hat{R}3)$, and this gives us the relation
$\tilde{h}_{21}(u)=\tilde{h}_{12}(u)^{-1}$. In addition we have $1=\tilde{c}(u,u^{-1})=\tilde{h}_{12}(u)\tilde{h}_{12}(u^{-1})$, that is,
$\tilde{h}_{12}(u)^{-1}=\tilde{h}_{12}(u^{-1})$ by $(ST3)$. Therefore we see $\tilde{z}$ is of the form
\[ \tilde{z}=\tilde{h}_{12}(u_1)\dotsm\tilde{h}_{12}(u_r)\]
for $u_i\in\Dt^{\times}$. We note that $\tilde{\phi}(\tilde{z})=I$ implies $u_1u_2\cdots u_r=1$ and $u_ru_{r-1}\cdots u_1=1$.
In addition, we easily see $[u_1,u_2][u_2u_1,u_3]\dotsm [u_{r-1}\dotsm u_1,u_r]=1$ using a similar way in Section 3. Then, by $(ST3)$,
\begin{align*}
\tilde{z}&=\tilde{h}_{12}(u_1)\tilde{h}_{12}(u_2)\dotsm\tilde{h}_{12}(u_r)\\
	&=\tilde{c}(u_1,u_2)\tilde{h}_{12}(u_2u_1)\tilde{h}_{12}(u_3)\dotsm\tilde{h}_{12}(u_r)\\
	&=\tilde{c}(u_1,u_2)\tilde{c}(u_2u_1,u_3)\tilde{h}_{12}(u_3u_2u_1)\tilde{h}_{12}(u_4)\dotsm\tilde{h}_{12}(u_r)\\
	&\hspace{20mm}\vdots\\
	&=\tilde{c}(u_1,u_2)\tilde{c}(u_2u_1,u_3)\dotsm\tilde{c}(u_{r-2}\dotsm u_1,u_{r-1})\tilde{h}_{12}(u_{r-1}\dotsm u_1)
\tilde{h}_{12}(u_r)\\
	&=\tilde{c}(u_1,u_2)\tilde{c}(u_2u_1,u_3)\dotsm\tilde{c}(u_{r-1}\dotsm u_1,u_r)\tilde{h}_{12}(u_r\dotsm u_1)\\
	&=\tilde{c}(u_1,u_2)\tilde{c}(u_2u_1,u_3)\dotsm\tilde{c}(u_{r-1}\dotsm u_1,u_r)\\
	&=1.\hspace{5mm}\square
\end{align*}\vspace{3mm}

\textbf{STEP 2.} We claim $\tilde{y}=1$.

We set subgroups of $U$ as follows:
\[ U_1=\{ e_{\alpha}(f)|f\in D[t]\}\ \mathrm{and}\ U_2=\{ e_{-\alpha}(g)|g\in D[t]t\}.\]
Then $U=\langle U_1,U_2\rangle$, but we will show that this is actually a free product.
We now introduce a ``degree map'' $\mathrm{deg}: D[t]\to \mathbb{Z}_{\geq 0}\cup\{ \infty\}$, which is defined by
$\mathrm{deg}(f)=m$ for $f=\Sigma_{i=0}^mf_it^i$ and $\mathrm{deg}(0)=\infty$, where $f_i\in D$.\vspace{3mm}

Let $x\in U$, and we assume
\[x=e_{\beta_1}(q_1)e_{\beta_2}(q_2)\dotsm e_{\beta_r}(q_r),\] 
where $\beta_i\in\Delta$, $\beta_i\neq\beta_{i+1}$, $q_i\neq 0$ and
\[q_i\in
		\begin{cases}
			D[t]\ \mathrm{if}\ \beta_i\in\Delta^+,\\
			D[t]t\ \mathrm{if}\ \beta_i\in\Delta^-.
		\end{cases}\]
For each $1\leq i\leq r$, we put
\[\left(\begin{array}{cc}
 a_i & b_i \\
 c_i & d_i
\end{array}\right)
=e_{\beta_1}(q_1)e_{\beta_2}(q_2)\dotsm e_{\beta_r}(q_r).\]
Suppose $\beta_1=\epsilon_1-\epsilon_2$. Then we have
\[\left(\begin{array}{cc}
 a_1 & b_1 \\
 c_1 & d_1
\end{array}\right)=
\left(\begin{array}{cc}
 1 & q_1 \\
 0 & 1
\end{array}\right)\]
and $0=\mathrm{deg}(a_1)\leq\mathrm{deg}(b_1)$. Since
\[\left(\begin{array}{cc}
 a_2 & b_2 \\
 c_2 & d_2
\end{array}\right)=
\left(\begin{array}{cc}
 a_1 & b_1 \\
 c_1 & d_1
\end{array}\right)
\left(\begin{array}{cc}
 1 & 0 \\
 q_2 & 1
\end{array}\right)=
\left(\begin{array}{cc}
 a_1+b_1q_2 & b_1 \\
 c_1+d_1q_2 & d_1
\end{array}\right),\]
we obtain $\mathrm{deg}(a_2)>\mathrm{deg}(b_2)=\mathrm{deg}(b_1)$. Moreover, since
\[\left(\begin{array}{cc}
 a_3 & b_3 \\
 c_3 & d_3
\end{array}\right)=
\left(\begin{array}{cc}
 a_2 & b_2 \\
 c_2 & d_2
\end{array}\right)
\left(\begin{array}{cc}
 1 & q_3 \\
 0 & 1
\end{array}\right)
\left(\begin{array}{cc}
 a_2 & a_2q_3+b_2 \\
 c_2 & c_2q_3+d_2
\end{array}\right),\]
we have $\mathrm{deg}(a_2)=\mathrm{deg}(a_3)\leq\mathrm{deg}(b_3)$. Continuing this process,
we can reach
\begin{align*}
&\mathrm{deg}(a_1)\leq\mathrm{deg}(b_1)=\mathrm{deg}(b_2)<\mathrm{deg}(a_2)=\\
&\mathrm{deg}(a_3)\leq\mathrm{deg}(b_3)=\mathrm{deg}(b_4)<\mathrm{deg}(a_4)=\\
&\mathrm{deg}(a_5)\leq\mathrm{deg}(b_5)=\mathrm{deg}(b_6)<\mathrm{deg}(a_6)=\dotsm .
\end{align*}
Next, suppose $\beta_1=\epsilon_2-\epsilon_1$. Then, in the same way as above, we can get
\begin{align*}
&\mathrm{deg}(b_1)\leq\mathrm{deg}(a_1)=\mathrm{deg}(a_2)<\mathrm{deg}(b_2)=\\
&\mathrm{deg}(b_3)\leq\mathrm{deg}(a_3)=\mathrm{deg}(a_4)<\mathrm{deg}(b_4)=\\
&\mathrm{deg}(b_5)\leq\mathrm{deg}(a_5)=\mathrm{deg}(a_6)<\mathrm{deg}(b_6)=\dotsm .
\end{align*}
In any case we can show $x\neq 1$, which means that $U=U_1\star U_2$ is a free product,
and we have $\tilde{\phi}(\tilde{y})=I$ implies $\tilde{y}=1$. $\square$\vspace{3mm}

{\prop{Notation is as above. We obtain that
\[ K_2(2,\Dt)=\langle \hat{c}(u_1,v_1)\hat{c}(u_2,v_2)\dotsm\hat{c}(u_r,v_r) | u_i,v_i\in\Dt^{\times};
[u_1,v_1][u_2,v_2]\dotsm [u_r,v_r]=1\rangle,\]
where $\hat{c}(u,v)=\hat{h}_{12}(u)\hat{h}_{12}(v)\hat{h}_{12}(vu)^{-1}$, and
$K_2(2,\Dt)$ is a central subgroup of $St(2,\Dt)$.}}\vspace{5mm}

\textbf{Proof:} It suffices to check that $\hat{c}(u_1,v_1)\hat{c}(u_2,v_2)\dotsm\hat{c}(u_r,v_r)$ is a
central element when $[u_1,v_1][u_2,v_2]\dotsm [u_r,v_r]=1$.
We easily see for $f\in\Dt$ and $1\leq i\neq j\leq 2$
\begin{align*}
\{ \hat{c}(u_1,v_1)\hat{c}(u_2,v_2)&\dotsm\hat{c}(u_r,v_r)\} \hat{x}_{ij}(f)
\{ \hat{c}(u_1,v_1)\hat{c}(u_2,v_2)\dotsm\hat{c}(u_r,v_r)\} ^{-1}\\
	&=\hat{x}_{ij}(\{ [u_1,v_1][u_2,v_2]\dotsm [u_r,v_r]\} ^{\pm 1}f\{ [u_1,v_1][u_2,v_2]\dotsm [u_r,v_r]\} ^{\mp 1})\\
	&=\hat{x}_{ij}(f).\ \square
\end{align*}\vspace{3mm}

Using the Proposition 4 and \cite{mj}, we construct a central extension of higher rank.
Suppose $n\geq 3$. We now consider the following commutative diagram;

\begin{align*}
1 \rightarrow K_2&(2,D_\tau) \rightarrow St(2,D_{\tau}) \rightarrow E(2,D_{\tau}) \rightarrow 1 \hspace{5mm}(\text{exact})\\
&\downarrow\hspace{20mm}\downarrow\hspace{20mm}\downarrow\\
1 \rightarrow K_2&(n,D_\tau) \rightarrow St(n,D_{\tau}) \rightarrow E(n,D_{\tau}) \rightarrow 1 \hspace{5mm}(\text{exact})
\end{align*}

Then we see that a canonical homomorphism of $K_2(2,\Dt)$ into $K_2(n,\Dt)$ is surjective by \cite{mj}
since $\Dt$ is a euclidean ring, so we have;\vspace{3mm}

{\thm{Let $\tilde{E}(n,\Dt)$ be the group presented by the generators $\tilde{x}_{ij}(f)$ for all $1\leq i\neq j\leq n$ and $f\in\Dt$ with the defining relations $(ST1)$, $(ST2)$ and $(ST3)$. Then $\tilde{E}(n,\Dt)\cong E(n,\Dt)$. $\square$}}\vspace{3mm}\\

{\prop{We have
\[ K_2(n,\Dt)=\langle \hat{c}(u_1,v_1)\hat{c}(u_2,v_2)\dotsm\hat{c}(u_r,v_r) | u_i,v_i\in\Dt^{\times};
[u_1,v_1][u_2,v_2]\dotsm [u_r,v_r]=1\rangle,\]
where $\hat{c}(u,v)=\hat{h}_{12}(u)\hat{h}_{12}(v)\hat{h}_{12}(vu)^{-1}$, and
$K_2(n,\Dt)$ is a central subgroup of $St(n,\Dt)$. $\square$}}\vspace{3mm}

\section{Universal central extensions}

We will show that the natural homomorphism $\phi$ is universal when $|Z(D)|\geq 5$ and $|Z(D)|\neq 9$, where
$Z(D)$ is the center of $D$.
We first fix four central elements $a,b,c,d\in D$, which satisfy
\[ a^2-1\neq 0,\ b=(a^2-1)^{-1},\ c\neq 0,\ c-1\neq 0,\ c^2-c+1\neq 0,\ d^3-1\neq 0.\]
We note that $E(n,\Dt)$ is perfect (also $St(n,\Dt)$) since
\[ x_{\beta}(f)=[h_{\beta}(a),x_{\beta}(bf)],\]
where $\beta\in\Delta$, $f\in\Dt$. Let $\phi^*:E^*\to E(n,\Dt)$ be any central extension, and set
\[M(z)=\{ z^*\in E^* | \phi^*(z^*)=z\} \]
for $z\in E(n,\Dt)$. Then we define
\[x_{\beta}^*(f)=[h^*,x^*]\]
for $h^*\in M(h_{\beta}(a))$ and $x^*\in M(x_{\beta}(bf))$. We see this is well-defined from
the following lemma, which is called ``central trick''.\vspace{3mm}

{\lem{Let $\psi: K\to G$ be a central extension. Then for $X,X',Y,Y'\in K $, we obtain that $\psi(X)=\psi(X')$ and $\psi(Y)=\psi(Y')$ imply $[X,Y]=[X',Y']$.}}\vspace{3mm}\\
\textbf{Proof:} This holds by direct calculation. $\square$\vspace{3mm}

For $u\in \Dt$ we put
\begin{align*}
w_{\beta}^*(u)&=x_{\beta}^*(u)x_{\beta}^*(-u^{-1})x_{\beta}^*(u),\\
h_{\beta}^*(u)&=w_{\beta}^*(u)w_{\beta}^*(-1).
\end{align*}

{\lem{Let $\beta\in\Delta$, $f\in\Dt$ and $u\in\Dt^{\times}$. Then the following equations hold.
\begin{align*}
&(1)\ w_{\beta}^*(u)x_{\beta}^*(f)w_{\beta}^*(u)^{-1}=x_{-\beta}^*(-u^{-1}fu^{-1})\\
&(2)\ h_{\beta}^*(u)x_{\beta}^*(f)h_{\beta}^*(u)^{-1}=x_{\beta}^*(ufu)\\
\end{align*}}}\vspace{3mm}
\textbf{Proof:}
\begin{align*}
(1)\ w_{\beta}^*(u)x_{\beta}^*(f)w_{\beta}^*(u)^{-1}
&=w_{\beta}^*(u)[h_{\beta}^*(a),x_{\beta}^*(bf)]w_{\beta}^*(u)^{-1}\\
&=[w_{\beta}^*(u)h_{\beta}^*(a)w_{\beta}^*(u)^{-1},w_{\beta}^*(u)x_{\beta}^*(bf)w_{\beta}^*(u)^{-1}]\\
&=[h_{-\beta}^*(a),x_{-\beta}^*(-bu^{-1}fu^{-1})]\ \text{(by central trick)}\\
&=x_{-\beta}^*(-u^{-1}fu^{-1}).
\end{align*}
\begin{align*}
(2)\ h_{\beta}^*(u)x_{\beta}^*(f)h_{\beta}^*(u)^{-1}
&=h_{\beta}^*(u)[h_{\beta}^*(a),x_{\beta}^*(bf)]h_{\beta}^*(u)^{-1}\\
&=[h_{\beta}^*(u)h_{\beta}^*(a)h_{\beta}^*(u)^{-1},h_{\beta}^*(u)x_{\beta}^*(bf)h_{\beta}^*(u)^{-1}]\\
&=[h_{\beta}^*(a),x_{\beta}^*(bufu)]\ \text{(by central trick)}\\
&=x_{\beta}^*(ufu).\hspace{10mm}\square
\end{align*}\vspace{3mm}

We now define $\pi_{\beta,\gamma}(f,g)$ for $\beta,\gamma\in\Delta$ and $f,g\in \Dt$ by
\begin{align*}
\pi_{\beta,\gamma}(f,g)=
	\begin{cases}
		[x_{\beta}^*(f),x_{\gamma}^*(g)]x_{\beta+\gamma}^*(fg)^{-1}\hspace{5mm}&\text{if $j=k$},\\
		[x_{\beta}^*(f),x_{\gamma}^*(g)]x_{\beta+\gamma}^*(-gf)^{-1}&\text{if $i=l$},\\
		[x_{\beta}^*(f),x_{\gamma}^*(g)]&\text{otherwise},
	\end{cases}
\end{align*}
where $\beta=\epsilon_i-\epsilon_j$, $\gamma=\epsilon_k-\epsilon_l$, $i\neq j$ and $k\neq l$.
We note $\pi_{\beta,\gamma}(f,g)$ is central in $E^*$. Then we can show the following (cf.\cite{s}).\vspace{3mm}

{\lem{Let $\beta=\epsilon_i-\epsilon_j,\gamma=\epsilon_k-\epsilon_l\in\Delta$ and $f,f',g,g'\in\Dt$. Then
\begin{align*}
&(1)\ \pi_{\beta,\gamma}(f+f',g)=\pi_{\beta,\gamma}(f,g)\pi_{\beta,\gamma}(f',g),\\
&(2)\ \pi_{\beta,\gamma}(f,g+g')=\pi_{\beta,\gamma}(f,g)\pi_{\beta,\gamma}(f,g'),\\
&(3)\ \pi_{\beta,\gamma}(f,g)=1,\\
&(4)\ x_{\beta}^*(f)x_{\beta}^*(g)=x_{\beta}^*(f+g).
\end{align*}}}
\textbf{Proof:} We will prove this lemma dividing into the three steps: we first check (1), (2) and (3) only for
$j\neq k$ and $i\neq l$, and then show that (4) holds for all $\beta\in\Delta$, $f,g\in\Dt$.
Finally, we check (1), (2) and (3) for the others.\vspace{3mm}

\textbf{Step 1.} Let $j\neq k,i\neq l$. Then, by definition of $\pi_{\beta,\gamma}(f,g)$,
\begin{align*}
\pi_{\beta,\gamma}(f+f',g)&=[x_{\beta}^*(f+f'),x_{\gamma}^*(g)]\\
		&=[x_{\beta}^*(f)x_{\beta}^*(f'),x_{\gamma}^*(g)]\\
		&=x_{\beta}^*(f)x_{\beta}^*(f')x_{\gamma}^*(g)x_{\beta}^*(f')^{-1}x_{\beta}^*(f)^{-1}x_{\gamma}^*(g)^{-1}\\
		&=x_{\beta}^*(f)\pi_{\beta,\gamma}(f',g)x_{\gamma}^*(g)x_{\beta}^*(f)^{-1}x_{\gamma}^*(g)^{-1}\\
		&=\pi_{\beta,\gamma}(f,g)\pi_{\beta,\gamma}(f',g).
\end{align*}
\begin{align*}
\pi_{\beta,\gamma}(f,g+g')&=[x_{\beta}^*(f),x_{\gamma}^*(g+g')]\\
		&=[x_{\beta}^*(f),x_{\gamma}^*(g)x_{\gamma}^*(g')]\\
		&=x_{\beta}^*(f)x_{\gamma}^*(g)x_{\gamma}^*(g')x_{\beta}^*(f)^{-1}x_{\gamma}^*(g')^{-1}x_{\gamma}^*(g)^{-1}\\
		&=x_{\beta}^*(f)x_{\gamma}^*(g)x_{\beta}^*(f)^{-1}\pi_{\beta,\gamma}(f,g')x_{\gamma}^*(g)^{-1}\\
		&=\pi_{\beta,\gamma}(f,g)\pi_{\beta,\gamma}(f,g').
\end{align*}\vspace{3mm}

If $i\neq k$ and $j\neq l$, then
\begin{align*}
\pi_{\beta,\gamma}(bf,g)&=h^*_{\beta}(a)\pi_{\beta,\gamma}(bf,g)h^*_{\beta}(a)^{-1}\\
	&=[h^*_{\beta}(a)x^*_{\beta}(bf)h^*_{\beta}(a)^{-1}, h^*_{\beta}(a)x^*_{\gamma}(g)h^*_{\beta}(a)^{-1}]\\
	&=[x^*_{\beta}(a^2bf),x^*_{\gamma}(g)]\hspace{5mm}\text{(by central trick and (R4))}\\
	&=\pi_{\beta,\gamma}(a^2bf,g)
\end{align*}
and
\begin{align*}
\pi_{\beta,\gamma}(f,g)&=\pi_{\beta,\gamma}(a^2bf-bf,g)\\
	&=[x^*_{\beta}(a^2bf)x^*_{\beta}(bf)^{-1},x^*_{\gamma}(g)]\\
	&=x^*_{\beta}(a^2bf)x^*_{\beta}(bf)^{-1}\pi_{\beta,\gamma}(bf,g)^{-1}x^*_{\beta}(bf)x^*_{\gamma}(g)x^*_{\beta}(a^2bf)^{-1}x^*_{\gamma}(g)^{-1}\\
	&=\pi_{\beta,\gamma}(a^2bf,g)\pi_{\beta,\gamma}(bf,g)^{-1},
\end{align*}
so we obtain $\pi_{\beta,\gamma}(f,g)=1$.\vspace{3mm}

If $i=k$ and $j\neq l$, then
\begin{align*}
\pi_{\beta,\gamma}(f,g)&=h^*_{\beta}(d)\pi_{\beta,\gamma}(f,g)h^*_{\beta}(d)^{-1}\\
	&=\pi_{\beta,\gamma}(d^2f,dg)\hspace{5mm}\text{(by central trick and (R4))}
\end{align*}
and
\begin{align*}
\pi_{\beta,\gamma}(f,g)&=h^*_{\beta-\gamma}(d)\pi_{\beta,\gamma}(f,g)h^*_{\beta-\gamma}(d)^{-1}\\
	&=\pi_{\beta,\gamma}(df,d^{-1}g)\hspace{5mm}\text{(by central trick and (R4))},
\end{align*}
so we obtain $\pi_{\beta,\gamma}(f,g)=\pi_{\beta,\gamma}(d^2f,dg)=\pi_{\beta,\gamma}(d^3f,g)$,
hence $\pi_{\beta,\gamma}((d^3-1)f,g)=1$. Therefore $\pi_{\beta,\gamma}(f,g)=1$.\vspace{3mm}

If $i\neq k$ and $j=l$, then using the same way as the previous case,
we see \[\pi_{\beta,\gamma}((d^3-1)f,g)=1.\] Hence $\pi_{\beta,\gamma}(f,g)=1$.\vspace{3mm}

If $i=k$ and $j=l$, then 
\[\pi_{\beta,\beta}(f,g)=h_{\beta}^*(c)\pi_{\beta,\beta}(f,g)h_{\beta}^*(c)^{-1}=[x_{\beta}^*(c^2f),x_{\beta}^*(c^2g)]=\pi_{\beta,\beta}(c^2f,c^2g),\]
so we see
\begin{align*}
\pi_{\beta,\gamma}(f,g)&=\pi_{\beta,\gamma}(cf+(1-c)f,g)\\
			&=\pi_{\beta,\gamma}(cf,g)\pi_{\beta,\gamma}((1-c)f,g)\\
			&=\pi_{\beta,\gamma}(f,g/c)\pi_{\beta,\gamma}(f,g/(1-c))\\
			&=\pi_{\beta,\gamma}(f,g/c+g/(1-c))\\
			&=\pi_{\beta,\gamma}(f,g/c(1-c))\\
			&=\pi_{\beta,\gamma}(c(1-c)f,g),
\end{align*}
hence $\pi_{\beta,\gamma}((c^2-c+1)f,g)=1$. Therefore we obtain $\pi_{\beta,\gamma}(f,g)=1$, in particular,
\[ x^*_{\beta}(f)x^*_{\beta}(g)=x^*_{\beta}(g)x^*_{\beta}(f)\]
for all $f,g\in\Dt$.\vspace{3mm}

\textbf{Step 2.} We now put $(f,g)=x_{\beta}^*(f)x_{\beta}^*(g)x_{\beta}^*(f+g)^{-1}$. Then, by Step 1,
\begin{align*}
(bf,bg)&=h_{\beta}^*(a)(bf,bg)h_{\beta}^*(a)^{-1}\\
	&=h_{\beta}^*(a)x_{\beta}^*(bf)x_{\beta}^*(bg)x_{\beta}^*(b(f+g))^{-1}h_{\beta}^*(a)^{-1}\\
	&=[h_{\beta}^*(a),x_{\beta}^*(bf)]x_{\beta}^*(bf)h_{\beta}^*(a)x_{\beta}^*(bg)x_{\beta}^*(b(f+g))^{-1}h_{\beta}^*(a)^{-1}\\
	&=[h_{\beta}^*(a),x_{\beta}^*(bf)]x_{\beta}^*(bf)[h_{\beta}^*(a),x_{\beta}^*(bg)]x_{\beta}^*(bg)\\
			&\hspace{10mm}\cdot [h_{\beta}^*(a),x_{\beta}^*(b(f+g))^{-1}]x_{\beta}^*(b(f+g))^{-1}\\
	&=x_{\beta}^*(f)x_{\beta}^*(bf)x_{\beta}^*(g)x_{\beta}^*(bg)x_{\beta}^*(f+g)^{-1}x_{\beta}^*(b(f+g))^{-1}\\
	&=(f,g)(bf,bg).
\end{align*}
Thus $(f,g)=1$ for all $f,g\in\Dt$, hence $x_{\beta}^*(f)x_{\beta}^*(g)=x_{\beta}^*(f+g)$.\vspace{3mm}

\textbf{Step 3.} Let $j=k$. Then, by Step 1 and Step 2,
\begin{align*}
\pi_{\beta,\gamma}(f+f',g)&=[x^*_{\beta}(f+f'),x^*_{\gamma}(g)]x^*_{\beta+\gamma}(-(f+f')g)\\
	&=x^*_{\beta}(f)x^*_{\beta}(f')x^*_{\gamma}(g)x^*_{\beta}(-f')x^*_{\beta}(-f)x^*_{\gamma}(-g)x^*_{\beta+\gamma}(-fg)x^*_{\beta+\gamma}(-f'g)\\
	&=x^*_{\beta}(f)[x^*_{\beta}(f'),x^*_{\gamma}(g)]x^*_{\gamma}(g)x^*_{\beta}(-f)x^*_{\gamma}(-g)x^*_{\beta+\gamma}(-fg)x^*_{\beta+\gamma}(-f'g)\\
	&=x^*_{\beta}(f)\pi_{\beta,\gamma}(f',g)x^*_{\beta+\gamma}(f'g)x^*_{\beta}(-f)\pi_{\beta,\gamma}(f,g)x^*_{\beta+\gamma}(-f'g)\\
	&=\pi_{\beta,\gamma}(f,g)\pi_{\beta,\gamma}(f',g)
\end{align*}
and we similarly obtain
\begin{align*}
\pi_{\beta,\gamma}(f,g+g')&=[x^*_{\beta}(f),x^*_{\gamma}(g+g')]x^*_{\beta+\gamma}(-f(g+g'))\\
	&=x^*_{\beta}(f)x^*_{\gamma}(g)x^*_{\gamma}(g')x^*_{\beta}(-f)x^*_{\gamma}(-g')x^*_{\gamma}(-g)x^*_{\beta+\gamma}(-fg)x^*_{\beta+\gamma}(-fg')\\
	&=x^*_{\beta}(f)x^*_{\gamma}(g')x^*_{\beta}(-f)[x^*_{\beta}(f),x^*_{\gamma}(g)]x^*_{\gamma}(-g')x^*_{\beta+\gamma}(-fg)x^*_{\beta+\gamma}(-fg')\\
	&=\pi_{\beta,\gamma}(f,g')x^*_{\beta+\gamma}(fg')\pi_{\beta,\gamma}(f,g)x^*_{\gamma}(-g)x^*_{\beta+\gamma}(-fg)x^*_{\beta+\gamma}(-fg')\\
	&=\pi_{\beta,\gamma}(f,g)\pi_{\beta,\gamma}(f,g').
\end{align*}

Also, we see
\begin{align*}
\pi_{\beta,\gamma}(f,g)&=h^*_{\beta}(d)\pi_{\beta,\gamma}(f,g)h^*_{\beta}(d)^{-1}\\
	&=h^*_{\beta}(d)[x^*_{\beta}(f),x^*_{\gamma}(g)]x^*_{\beta+\gamma}(-fg)h^*_{\beta}(d)^{-1}\\
	&=[x^*_{\beta}(d^2f),x^*_{\gamma}(d^{-1}g)]x^*_{\beta+\gamma}(-dfg)\\
	&=\pi_{\beta,\gamma}(d^2f,d^{-1}g)
\end{align*}
and
\begin{align*}
\pi_{\beta,\gamma}(f,g)&=h^*_{\beta+\gamma}(d)\pi_{\beta,\gamma}(f,g)h^*_{\beta+\gamma}(d)^{-1}\\
	&=[x^*_{\beta}(df),x^*_{\gamma}(dg)]x^*_{\beta+\gamma}(-d^2fg)\\
	&=\pi_{\beta,\gamma}(df,dg).
\end{align*}
Therefore we obtain $\pi_{\beta,\gamma}((d^3-1)f,g)=1$, in particular, $\pi_{\beta,\gamma}(f,g)=1$.\vspace{3mm}

Let $i=l$. Then, using the similar way to the case $j=k$, we obtain
\begin{align*}
\pi_{\beta,\gamma}(f+f',g)&=\pi_{\beta,\gamma}(f,g)\pi_{\beta,\gamma}(f',g),\\
\pi_{\beta,\gamma}(f,g+g')&=\pi_{\beta,\gamma}(f,g)\pi_{\beta,\gamma}(f,g'),\\
\pi_{\beta,\gamma}(f,g)&=1.\ \square
\end{align*}\vspace{3mm}


From the above, we have the following (cf. [7, p.51]).

{\thm{If $|Z(D)|\geq 5$ and $|Z(D)|\neq 9$, then $\phi$ is a universal central extension of $E(n,\Dt)$.}}\vspace{3mm}

\textbf{Proof:} For any central extension $\phi^*$, we can construct $x_{\beta}^*(f)$ for all $\beta\in\Phi$ and $f\in\Dt$
as above. Then the relations give us a homomorphism $\hat{\phi}^*:St(n,\Dt)\to E^*$ defined by
$\hat{\phi}^*(e_{\beta}(f))=x_{\beta}^*(f)$. Thus we obtain $\phi=\phi^*\circ \hat{\phi}^*$. $\square$\vspace{5mm} 

\begin{rem}
If the cardinality of $Z(D)$ is less than 5 or is 9, then there exist counter examples which are written in \cite{s}.
\end{rem}

\begin{rem}
Let $F$ be any field and $x_i$ $(i\in\mathbb{Z})$ be countable indeterminates. We put
\[K=F(x_i)_{i\in\mathbb{Z}}=\left\{ \frac{f(x_{i_1},\cdots,x_{i_r})}{g(x_{j_1},\dots,x_{j_s})}\ \Bigg|\ f,g\in F[x_i]_{i\in\mathbb{Z}} \right\}.\]

We define an automorphism $\tau$ of $K$ to be $\tau(x_i)=x_{i+1}$. Furthermore let $t$ be a new indeterminate, and put
\[ D=K((t))=\left\{ \sum_{-\infty<<k<\infty}a_kt^k\ \Bigg|\ a_k\in K\right\}. \] 
Then, $D$ can be viewed as a division ring, whose product is given by
\[\left(\sum_ka_kt^k\right)\left(\sum_lb_lt^l\right)=\sum_m\left(\sum_{k+l=m}a_k\tau^k(b_l)\right)t^m.\]
In particular,  for any field $F$ there exists a division ring $D$ such that $Z(D)=F$.
\end{rem}

\section{$K_1$-groups}
As a subgroup of $GL(n,\Dt)$, if we put $H=\{ \mathrm{diag}(u_1,\dots, u_n)| u_i\in \Dt^{\times}\}$, then
$GL(n,\Dt)=\langle E(n,\Dt),H\rangle$ since $\Dt$ is a euclidean ring (cf. \cite{ms},\cite{hs} Proposition 1.1.2).
Put $H_1=\{ \mathrm{diag}(u,1,\dots,1)|u\in\Dt^{\times} \}$, and we have
\begin{align*}
GL(n,\Dt)&=\langle E(n,\Dt), H\rangle\\
		&=\langle E(n,\Dt), H_1\rangle\\
		&\rhd E(n,\Dt).
\end{align*}
Define our $K_1$-group by $K_1(n,\Dt)=GL(n,\Dt)/E(n,\Dt)$ (cf. \cite{m}). Then, we have
\[ K_1(n,\Dt)=GL(n,\Dt)/E(n,\Dt)\cong H_1/(E(n,\Dt)\cap H_1).\]
In the following we will discuss the subgroup $E(n,\Dt)\cap H_1$.

Since $E(n,\Dt)$ has a Bruhat decomposition, we consider $BwB\cap H_1$. From Proposition 1,
there exist $\dot{w}\in W$ and $h=\mathrm{diag}(t^{m_1},\dots,t^{m_n})$ such that $w=\dot{w}h$,
where $m_i\in\mathbb{Z}$ with $m_1+\cdots+m_n=0$, and then
\[BwB\cap H_1=UT_0wT_0U\cap H_1=U\dot{w}hT_0U\cap H_1.\]
If we set
\[\Dt^{\geq 0}=D[t]\ \text{and}\ \Dt^{>0}=D[t]t,\]
then we see
\[U\subset
\left(
\begin{array}{cccc}
1+\Dt^{>0} & \Dt^{\geq 0} & \cdots & \Dt^{\geq 0}\\
\Dt^{>0}    & 1+\Dt^{>0}    &           & \vdots\\
\vdots       &                   & \ddots & \Dt^{\geq 0}\\
\Dt^{>0}    & \cdots         & \Dt^{>0} &1+\Dt^{>0}
\end{array}
\right).\]
Suppose $BwB\cap H_1\neq\emptyset$. Then
\begin{align*}
BwB\cap H_1\neq\emptyset
&\Rightarrow U\dot{w}hT_0U\cap H_1\neq\emptyset\\
&\Rightarrow U\cap H_1UT_0h^{-1}{\dot{w}}^{-1}\neq\emptyset.
\end{align*}
Therefore we can take $d\in H_1$ of degree $m$, $x_+\in U$ and
$h_0=\mathrm{diag}(u_1,\dots,u_n)\in T_0$ satisfying
$x=dx_+h_0h^{-1}{\dot{w}}^{-1}\in U$. Then we see
\[x\in
\left(
\begin{array}{cccc}
t^{m-m_1}(1+\Dt^{>0}) & t^{m-m_2}\Dt^{\geq 0} & \cdots & t^{m-m_n}\Dt^{\geq 0} \\
t^{-m_1}\Dt^{>0} & t^{-m_2}(1+\Dt^{>0}) & \cdots & t^{-m_n}\Dt^{\geq 0}\\
\vdots & \vdots & \ddots & \vdots\\
t^{-m_1}\Dt^{>0} & t^{-m_2}\Dt^{>0} & \cdots & t^{-m_n}(1+\Dt^{>0})
\end{array}
\right)
\dot{w}^{-1}.\]

If $m_i>0$ for some $2\leq i\leq n$, then $x$ has an entry with a negative power of $t$ in a diagonal part.
This is a contradiction, hence $m_i\leq 0$ for all $2\leq i\leq n$. In particular we obtain $m_1\geq 0$ from $m_1+\cdots+m_n=0$.
Suppose that $\dot{w}^{-1}$ is a permutation sending $j$th column to the 1st one for some $2\leq j\leq n$.
Then $2\leq \dot{w}^{-1}(k)\leq n$ for all $2\leq k\leq n$ without $k= j$.
If $\dot{w}^{-1}(i_1)>\dot{w}^{-1}(i_2)$ for $2\leq i_1<i_2\leq n$ except $j$,
then we see $t^{-m_{i_1}}\Dt^{>0}\cap \Dt^{\geq 0}=\emptyset$. This is a contradiction.
Thus, we obtain $\dot{w}^{-1}(k)=k$ for such $k$.

Hence, it suffices to check the two cases: When $\dot{w}^{-1}$ is a transposition of the 1st column and the $j$th column,
and when $\dot{w}^{-1}$ is an identity. Suppose that $\dot{w}^{-1}$ is a transposition $(1,j)$.
Then we have $t^{m-m_1}(1+\Dt^{>0})\cap\Dt^{\geq 0}\neq\emptyset$, and which implies $n-m_1\geq 0$.
On the other hand, we see $t^{m-m_j}\Dt^{\geq 0}\cap(1+\Dt^{>0})\neq\emptyset$, and which implies $m-m_j\leq 0$.
If we consider $m_1+m_j=0$, we obtain
\[ 0\leq -m_j\leq m\leq m_j\leq 0,\]
hence, $m=m_1=m_j=0$. However, the 1st column cannot be going to the $j$th column from the fact
$(1+\Dt^{>0})\cap\Dt^{>0}=\emptyset$. Therefore $\dot{w}=1$.

Thus, we obtain
\begin{align*}
E(n,\Dt)\cap H_1=B\cap H_1&=(U\rtimes T_0) \cap H_1\\
	&=T_0\cap H_1\\
	&\subset T\cap H_1.
\end{align*}

In the following, for $\beta=\epsilon_i-\epsilon_j$ we denote the element $h_{\beta}(s)$ in $T$ by $h_{ij}(s)$.
From the definition of $T$ and the relation $h_{ij}(s)=h_{1j}(s)h_{1i}(s^{-1})$, we find that $T$ is generated by
$h_{1j}(s)$ for all $2\leq j\leq n$ and $s\in\Dt^{\times}$.
Thus every element $h$ in $T\cap H_1$ can be of the form
\[ h=h_{1,l_1}(s_1)h_{1,l_2}(s_2)\cdots h_{1,l_k}(s_k),\]
where $2\leq l_1,\dots, l_k\leq n$ and $s_i\in\Dt^{\times}$.
Here, we put $I_j=\{ i| l_i=j\}$ for $2\leq j\leq n$. Then we see $\{ l_1,\dots, l_k\}=I_2\cup I_3\cup\cdots\cup I_n$.
Using this notation we find
\[h=\left(
\begin{array}{cccc}
\prod_{i=1}^ks_i & 0 & \cdots & 0 \\
0 & \prod_{i\in I_2}s_i^{-1} & & \huge{0}\\
\vdots & & \ddots & \\
0 & \huge{0} & & \prod_{i\in I_n}s_i^{-1}
\end{array}\right),\]
where $\prod_{i\in I_j}s_i^{-1}=1$ for all $2\leq j\leq n$.
Therefore, if we put $s=\prod_{i=1}^ks_i$, we can rewrite $s$ as
\[ s=s_1\cdots s_k=s_1\cdots s_ks_{i_1}^{-1}\dots s_{i_k}^{-1},\]
where $\{ i_1\dots,i_k\}=\{ 2,\dots,n\}$.
If $i_r=1$, then putting $v_1=s_2\cdots s_ks_{i_1}^{-1}\cdots s_{i_{r-1}}^{-1}$ we have
\begin{align*}
s&=[s_1,v_1]v_1s_1s_{i_r}^{-1}\cdots s_{i_k}^{-1}\\
	&=[s_1,v_1]v_1s_{i_{r+1}}^{-1}\cdots s_{i_k}^{-1}\\
	&=[s_1,v_1]s_2\cdots s_ks_{i_1}^{-1}\cdots s_{i_{r-1}}^{-1}s_{i_{r+1}}^{-1}\cdots s_{i_k}^{-1}.
\end{align*}
If $i_{r'}=2$, then putting $v_2=s_3\cdots s_ks_{i_1}^{-1}\cdots s_{i_{r'-1}}^{-1}$ we also have
\[ s=[s_1,v_1][s_2,v_2]s_3\cdots s_ks_{i_1}^{-1}\cdots s_{i_{r'-1}}^{-1}s_{i_{r'+1}}^{-1}\cdots s_{i_{r-1}}^{-1}s_{i_{r+1}}^{-1}\cdots s_{i_k}^{-1}.\]
Repeating this operation, we finally reach
\[ s=[s_1,v_1][s_2,v_2]\cdots [s_k,v_k].\]
Thus, we obtain
\[T\cap H_1\subset \left(
\begin{array}{cccc}
[\Dt^{\times},\Dt^{\times}] & 0 & \cdots & 0 \\
0 & 1 & & 0\\
\vdots & & \ddots & \\
0 & 0 & & 1
\end{array}\right).\]
On the other hand, we see
\[\left(
\begin{array}{cccc}
[u,v] & 0 & \cdots & 0 \\
0 & 1 & & 0\\
\vdots & & \ddots & \\
0 & 0 & & 1
\end{array}\right)=h_{12}(u)h_{12}(v)h_{12}(u^{-1}v^{-1}),\]
where $u,v\in\Dt^{\times}$. We see $\mathrm{deg}([u,v])=0$, hence, 
\[ \left(
\begin{array}{cccc}
[\Dt^{\times},\Dt^{\times}] & 0 & \cdots & 0 \\
0 & 1 & & 0\\
\vdots & & \ddots & \\
0 & 0 & & 1
\end{array}\right) \subset T_0\cap H_1.\]
These deduce
\[E(n,\Dt)\cap H_1= \left(
\begin{array}{cccc}
[\Dt^{\times},\Dt^{\times}] & 0 & \cdots & 0 \\
0 & 1 & & 0\\
\vdots & & \ddots & \\
0 & 0 & & 1
\end{array}\right).\]
Therefore we obtain the following theorem.\vspace{3mm}

{\thm{Notation is as above. Then, we have
\[ K_1(n,\Dt)\cong \Dt^{\times}/[\Dt^{\times},\Dt^{\times}].\hspace{3mm}\square\]}}

\end{document}